\documentclass[psamsfonts,reqno,11pt]{amsart}
\usepackage{graphicx}
\usepackage{graphics}
\usepackage{latexsym}
\usepackage{amsmath}
\usepackage{amssymb}
\usepackage[numbers,sort&compress]{natbib}
\makeatletter
\@namedef{subjclassname@2010}{%
  \textup{2010} Mathematics Subject Classification}
\makeatother

\newtheorem{thm}{Theorem}[section]
\newtheorem{lem}[thm]{Lemma}
\newtheorem{eg}[thm]{Example}
\newtheorem{prop}[thm]{Proposition}

\newtheorem{defn}[thm]{Definition}
\newtheorem{rem}[thm]{Remark}

\newtheorem{rem-eg}[thm]{Remark and Example}

\newtheorem{problem}[thm]{Problem}

\newenvironment{prf}{{\noindent \textbf{Proof:} }}{\hfill $\Box$\medskip}
\def\sideremark#1{\ifvmode\leavevmode\fi\vadjust{\vbox
to0pt{\vss \hbox to 0pt{\hskip\hsize\hskip1em
\vbox{\hsize2cm\tiny\raggedright\pretolerance10000
\noindent#1\hfill}\hss}\vbox to8pt{\vfil}\vss}}}

\begin{document}

\title{Mappings of preserving $n$-distance one in $n$-normed spaces}

\author[X. Huang]{Xujian Huang}
\address{Department of Mathematics\\ Tianjin University of Technology\\ 300384 Tianjin\\ China}
\email{huangxujian86@sina.cn}

\author[D. Tan]{Dongni Tan$^{*}$}
\address{Department of Mathematics\\ Tianjin University of Technology\\ 300384 Tianjin\\ China}
\email{tandongni0608@sina.cn}

\thanks{
$*$Corresponding author.\\
\indent The authors are supported by the Natural Science Foundation of China (Grant Nos. 11201337, 11201338, 11371201, 11301384). The first author was supported by the Tianjin Science \& Technology Fund 20111001.}
\begin{abstract}
We give a positive answer to the Aleksandrov problem in $n$-normed spaces under the surjectivity assumption. Namely, we show that every surjective mapping preserving $n$-distance one is affine, and thus is an $n$-isometry. This is the first time to solve the Aleksandrov problem in $n$-normed spaces with only surjective assumption even in the usual case $n=2$. Finally, when the target space is $n$-strictly convex, we prove that every mapping preserving two $n$-distances with an integer ratio is an affine $n$-isometry.

\end{abstract}



\maketitle

\section{Introduction and Preliminaries}
\medskip

Let $X$ and $Y$ be two  metric spaces. A mapping $f: X \rightarrow Y$ is called an isometry if $f$ satisfies
\begin{eqnarray*}
d_{Y}(f(x),f(y))=d_{X}(x,y)
\end{eqnarray*}
for all $x,y\in X$, where $d_{X}(,)$ and $d_{Y}(,)$ denote the metric in the space $X$ and $Y$, respectively. For some $r>0$, suppose that $f$ preserves distance $r$, i.e., for all $x,y\in X$ with $d_{X}(x,y)=r$, we have $d_{Y}(f(x),f(y))=r$. Then $r$ is called a \emph{conservative distance} for the mapping $f$.

In 1970, Aleksandrov \cite{A} posed the following problem:
\begin{problem}
Under what conditions is a mapping of a metric space $X$ into itself preserving distance one an isometry?
\end{problem}
It is called the  Aleksandrov problem. It has been extensively investigated by many authors (see
\cite{CKK,CLP,G2,R1,R2,R3, R4, RS} and the references
therein).
This problem still remains open even in the case where $X = \mathbb{R}^n$ and $Y =\mathbb{R}^m$ with $2 < n < m$ (see \cite{R4}).

The study of $n$-normed spaces began early in the second half of the twentieth century (see \cite{G1,G,M1,M2}), and it is also an widely-studied and interesting area even today (see e.g. \cite{C,HPP,CKK,CLP}). Chu et al. \cite{CLP} first generalized the Aleksandrov problem to $n$-normed spaces. Their main result \cite[Theorem 2.10]{CLP} proves that the weak $n$-distance one preserving mapping is an $n$-isometries under additional conditions (e.g. $n$-1-Lipschitz, preserving 2-collinearity).

A natural question can be raised as a modified version of the Aleksandrov problem: What happens if two (or more) distances are preserved by a mapping between normed spaces? W. Benz \cite{B} (see also \cite{BB}) investigated the case when the mapping preserves distances $\rho$ and $n\rho$ for some $\rho>0$ and some integer $n>1$. If the target space is strictly convex, they showed in \cite{B} that this mapping is an affine isometry. If the mapping $f$ preserves two distances with a non-integer ratio, it is an open problem whether or not $f$ must be an isometry. For more information we refer to  \cite{R1,R2,R3,R4}. Motivated by these results and also as an application of our main results we shall show that the result of W. Benz remains valid in $n$-normed spaces if the target space is $n$-strictly convex.

In this paper, we show that every mapping between two $n$-normed spaces preserving a fixed nonzero weak $n$-distance and 2-collinearity for the midpoint of a segment is affine, and thus is an $n$-isometry. By this we show that every surjective mapping preserving $n$-distance one is an affine $n$-isometry.  Finally, if the target space is $n$-strictly convex, we show that every mapping preserves two $n$-distances with an integer ratio is an affine $n$-isometry.

Throughout this paper, all linear spaces will be assumed real. Let $n\geq 2$, $X$ and $Y$ be two linear $n$-normed spaces whose dimensions greater than $n-1$.

In the remainder of this introduction,  we will present some definitions in $n$-normed spaces and cite an example of $n$-normed spaces for the easy understanding of this kind of spaces.

An $n$-norm on a real vector space $X$ (of dimension at least $n$) is a mapping $\| \cdot, \cdots ,\cdot \|: X^{n} \rightarrow \mathbb{R}$ which satisfies the following four conditions:\\
{\rm(a)} $\| x_1, \cdots, x_n\|= 0$ if and only if $x_1, \cdots, x_n$ are linearly dependent;\\
{\rm(b)} $\| x_1, \cdots, x_n\|$ is invariant under permutation;\\
{\rm(c)} $\|\alpha x_1, \cdots, x_n\| = |\alpha| \| x_1, \cdots, x_n\|$ for $\alpha \in \mathbb{R}$;\\
{\rm(d)} $\| x_0+x_1, x_2, \cdots, x_n\| \leq  \| x_0,x_2, \cdots, x_n\| + \| x_1,x_2, \cdots, x_n\|$.\\
The pair $(X, \| \cdot, \cdots ,\cdot \|)$ is called an \emph{$n$-normed  space}. Note that in this space, we have $\|x_1, x_2, \cdots, x_n\|=\|x_1+y, x_2, \cdots, x_n\|$ for any linear combination $y$ of $x_2, \cdots, x_n\in X$.

\begin{eg}
If $X$ is a normed space with dual $X'$, then as formulated by \"{G}ahler \emph{(see \cite{G})}
we may define an $n$-norm on $X$ by
\begin{eqnarray*}
\| x_1, x_2, \cdots, x_n\| := \sup_{f_j\in X', \ \|f_j\|\leq 1}\begin{vmatrix}
f_1(x_1) & \cdots  & f_1(x_n)\\
\vdots   & \ddots  & \vdots\\
f_n(x_1) & \cdots  & f_n(x_n)\\
\end{vmatrix}=\sup_{f_j\in X', \ \|f_j\|\leq 1}\det[f_j(x_i)].
\end{eqnarray*}
Meanwhile, if $X$ is equipped with an inner product $\langle \cdot, \cdot \rangle$, we can define the standard
$n$-norm on $X$ by
\begin{eqnarray*}
\| x_1, x_2, \cdots, x_n\|:= \sqrt{det[\langle x_i, x_j\rangle]},
\end{eqnarray*}
which can be interpreted as the volume of the $n$-dimensional parallelepiped spanned by $x_1, x_2, \cdots, x_n\in X$ \emph{(see \cite{GH}).}
\end{eg}

Recall some definitions in $n$-normed spaces.

\begin{defn}
Let $X$ and $Y$ be two $n$-normed spaces, and let $f: X\rightarrow Y$ be a mapping.

\rm{(a)}  $f$ is said to be an \emph{$n$-isometry} if it satisfies
\begin{eqnarray*}
\| f(x_1)-f(y_1), \cdots, f(x_n)-f(y_n) \| = \| x_1-y_1, \cdots, x_n-y_n\|
\end{eqnarray*}
for all $x_1, \cdots, x_n, y_1, \cdots, y_n \in X$. In particular, if $ y_1=\cdots = y_n $, $f$ is said to be a \emph{weak $n$-isometry}.

\rm{(b)} $f$ is said to have the \emph{$n$-distance one preserving property} ($n$-DOPP), if
\begin{eqnarray*}\| x_1-y_1, \cdots, x_n-y_n\|=1 \Rightarrow \ \| f(x_1)-f(y_1), \cdots, f(x_n)-f(y_n) \| =1.\end{eqnarray*}
 for all $x_1, \cdots, x_n, y_1, \cdots, y_n \in X$. In particular, if $ y_1=\cdots = y_n $, $f$ is said to have the \emph{weak n-distance one preserving property}($w$-$n$-DOPP).

\rm{(c)} $f$ is said to \emph{preserve $\rho$-$n$-distance} for some $\rho>0$, if $\| x_1-y_1, \cdots, x_n-y_n\|=\rho$ implies $\| f(x_1)-f(y_1), \cdots, f(x_n)-f(y_n) \| =\rho$ for all $x_1, \cdots, x_n, y_1, \cdots, y_n \in X$.
In particular, if $ y_1=\cdots = y_n $, $f$ is said to \emph{preserve $w$-$\rho$-$n$-distance}.

\rm{(d)} $f$ is called an \emph{$n$-Lipschitz mapping} if there is a $K \geq 0$ such that
\begin{eqnarray*}
\| f(x_1)-f(y_1), \cdots, f(x_n)-f(y_n) \| \leq K \| x_1-y_1, \cdots, x_n-y_n\|
\end{eqnarray*}
for all $x_1, \cdots, x_n, y_1, \cdots, y_n \in X$. In this case, the constant $K$ is called the $n$-Lipschitz constant.
In particular, if $ y_1=\cdots = y_n $, $f$ is said to be a \emph{weak $n$-Lipschitz mapping}.
\end{defn}

\section{Isometry in $n$-normed spaces}

In this section we consider the Aleksandrov problem in $n$-normed spaces. We first introduce a weak case of preserving 2-collinearity.  Then, we prove that the Aleksandrov problem holds in $n$-normed spaces under weaker hypothesis.

Note that the points $x,y,z$ of $X$ are said to be \emph{2-collinear} if $y-z = t(x-z)$ for some real number $t$. The points $x_0,x_1,\cdots,x_n$ of $X$ are said to be \emph{$n$-collinear} if for some $i$, the points $x_j-x_i, 0\leq j\neq i\leq n$ are linearly dependent.

\begin{defn}
Let $X$ and $Y$ be two $n$-normed spaces, and let $f$ be a mapping from $X$ into $Y$.

\rm{(a)} $f$ is said to preserve \emph{2-collinearity} if $x,y,z\in X$ are collinear, then $ f(x),f(y),f(z)$ are collinear. In particular, if $z=(x + y)/2$, $f$ is said to preserve \emph{2-collinearity for the midpoint of a segment}.

{\rm(b)} $f$ is said to preserve \emph{$n$-collinearity} if $x_0, x_1,\cdots,x_n$ of $X$ are $n$-collinear, then $f(x_0), f(x_1),\\ \cdots, f(x_n)$ are $n$-collinear. That means that $f$ preserves $w$-$0$-distance, i.e., if $\|x_1-x_0,  \cdots, x_n-x_0\|=0$, then
\begin{eqnarray*}
\|f(x_1)-f(x_0),  \cdots, f(x_n)-f(x_0)\|=0\end{eqnarray*} for all $x_0, x_1, \cdots , x_n \in X$.
\end{defn}

In the first step, we prove the following lemma indicating that a mapping $f$ from an $n$-normed space $X$ to an $n$-normed space $Y$, which preserves a nonzero weak $n$-distance and 2-collinearity for the midpoint of a segment, satisfies Jensen's equation:
\begin{eqnarray*}
f(\frac{x+y}{2})= \frac{f(x)+f(y)}{2}, \quad \forall x,y \in X.
\end{eqnarray*}

\begin{lem}\label{add}
Let $X$ and $Y$ be two  $n$-normed spaces, and let $f : X \rightarrow Y$ preserve $w$-$\rho$-$n$-distance for some $\rho>0$. Then $f$ is injective.  Moreover if $f$ preserves 2-collinearity for the midpoint of a segment, then $f(x)-f(0)$ is additive.
\end{lem}
\begin{prf} For $x\neq y\in X$, the assumption that $\dim X\geq n$ allows the existence of $x_2, x_3, \cdots, x_n \in X$ such that
\begin{eqnarray*} \|y-x, x_2-x, \cdots, x_n-x\|=\rho. \end{eqnarray*}
Since the mapping $f$ preserves $w$-$\rho$-$n$-distance, we have
\begin{eqnarray*}
\| f(y)-f(x), f(x_2)-f(x),  \cdots, f(x_n)-f(x) \| = \rho.
\end{eqnarray*}
This implies $f(x) \neq f(y)$, and thus $f$ is injective. To see our second conclusion,
it suffices to prove that for all $x,y \in X$, we have
\begin{eqnarray}\label{eq1}
f(\frac{x+y}{2})= \frac{f(x)+f(y)}{2}.
\end{eqnarray}
To prove (\ref{eq1}), set $z = (x + y)/2$ for distinct $x,y\in X$.  Choose $y_2, y_3, \cdots, y_n \in X$ such that
\begin{eqnarray*} \|y-z, y_2-z, \cdots, y_n-z\|=\|x-z, y_2-z, \cdots, y_n-z\|=\rho. \end{eqnarray*}
Then clearly
\begin{eqnarray} \label{c1}
&&\|f(y)-f(z), f(y_2)-f(z),\cdots, f(y_n)-f(z)\|=\rho \\ \label{c2} &&\|f(x)-f(z), f(y_2)-f(z), \cdots, f(y_n)-f(z)\|=\rho.
\end{eqnarray}
Since $f$ preserves 2-collinearity for the midpoint of a segment, there exists a real number $t$ such that
\begin{eqnarray*}f(y)-f(z) = t (f(x) - f(z)).\end{eqnarray*}
By(\ref{c1}) and (\ref{c2}), we obtain that $t=-1$, and hence
\begin{eqnarray*}
f(\frac{x+y}{2})=f(z)= \frac{f(x)+f(y)}{2}.
\end{eqnarray*}
\end{prf}

One may wonder how to check that a mapping $f$ from an $n$-normed space into another preserves 2-collinearity.  What interests us is that it only requires $f$ to preserve $w$-$n$-DOPP (not necessarily surjective) and be a weak $n$-Lipschitz mapping or preserve $n$-collinearity. This has been indicated in \cite[Lemma 3.2]{CKK} which states that every $n$-isometry $f$ preserves 2-collinearity in $n$-normed spaces. For the convenience of readers and since the condition is weaker, we here include a proof.

\begin{lem}\label{col}
Let $X$ and $Y$ be two  $n$-normed spaces. Suppose that the mapping $f : X \rightarrow Y$ preserves $w$-$\rho$-$n$-distance for some $\rho>0$. Then the following properties are equivalent:\\
{\rm(a)} $f$ preserves $n$-collinearity;\\
{\rm(b)} $f$ preserves 2-collinearity;\\
{\rm(c)} $f$ preserves 2-collinearity for the midpoint of a segment.
\end{lem}

\begin{prf}
For the implication $(a)\Rightarrow (b)$ assume that, on the contrary, there are $x_0,x_1,x_2\in X$ which are collinear such that $f(x_1)-f(x_0), f(x_2)-f(x_0)$ are linearly independent. Note that $x_0\neq x_1$ and $f$ preserves $w$-$\rho$-$n$-distance. We can choose $y_2, \cdots, y_n \in X$ such that \begin{eqnarray*}&& \| f(x_1)-f(x_0), f(y_2)-f(x_0),  \cdots, f(y_n)-f(x_0) \|\\
 &&= \| x_1-x_0, y_2-x_0, \cdots, y_n-x_0\|=\rho.\end{eqnarray*}
Then the set $A:=\{f(x)-f(x_0): x\in X\}$ contains $n$ linearly independent vectors. Hence there exist $x_3, \cdots, x_n \in X$ such that
$$\|f(x_1)-f(x_0), f(x_2)-f(x_0), f(x_3)-f(x_0), \cdots, f(x_n)-f(x_0)\|\neq 0.$$  Assume that $f$ preserves $n$-collinearity. Then
$\|x_1-x_0, x_2-x_0,  \cdots, x_n-x_0\|=0$ implies that
\begin{equation*}
\|f(x_1)-f(x_0), f(x_2)-f(x_0), \cdots, f(x_n)-f(x_0)\|=0.
\end{equation*}
which is a contradiction. Thus $f$ preserves 2-collinearity.

The implication $(b)\Rightarrow (c)$ is clear.

For the implication $(c)\Rightarrow (a)$ without loss of generality we can assume that $\rho=1$. Then $f$ satisfies $w$-$n$-DOPP. Let $g(x)=f(x)-f(0)$ for every $x\in X$. We first prove that $g$ preserves distance $m/k$ for all $m,k \in \mathbb{N}$. Let $x_1, x_2, \cdots, x_n
$ be in $X$ and $m,k$  be in $\mathbb{N}$ such that $$\| x_1, x_2, \cdots, x_n\|=m/k.$$
We see from Lemma \ref{add} that $g$ is $\mathbb{Q}$-linear, and since $g(0)=0$ and satisfies $w$-$n$-DOPP, we have
\begin{eqnarray*}
\|g(x_1), g(x_2), \cdots, g(x_n)\|= \frac{m}{k}\|g(\frac{k}{m}x_1), g(x_2), \cdots, g(x_n)\|=\frac{m}{k}.
\end{eqnarray*}
 To see that $g$ preserves $n$-collinearity, we only need to check that for all $x_1, x_2, \cdots x_n\in X$ which are not all zero with $\|x_1, x_2, \cdots, x_n\|= 0$,
\begin{eqnarray*}
\| g(x_1), g(x_2), \cdots, g(x_n)\|= 0.
\end{eqnarray*}
Since $\|x_1, x_2, \cdots, x_n\|= 0$, we know that $x_1, x_2, \cdots, x_n\in X$ are linearly dependent.  To simplify the notation, the maximal linearly independent members of $x_1, x_2, \cdots, x_n$ are still denoted by $x_1,\cdots,x_k$ where $1\leq k<n$. Choose $y_{k+1},\cdots,y_n\in X$ such that
\begin{eqnarray*}\|x_1,\cdots,x_k, y_{k+1},\cdots,y_n\|=1.\end{eqnarray*}
Then for every positive integer $m$, \begin{eqnarray*}\|x_1,\cdots,x_k, x_{k+1}+\frac{1}{m}y_{k+1},\cdots, x_n+\frac{1}{m}y_n\|=\frac{1}{m^{n-k}}\end{eqnarray*} and by the above,
\begin{eqnarray*}
&&\| g(x_1),\cdots, g(x_k), g(x_{k+1})+\frac{1}{m}g(y_{k+1}),\cdots, g(x_n)+\frac{1}{m}g(y_{n})\|=\frac{1}{m^{n-k}}.
\end{eqnarray*}
Triangle inequality hence gives
\begin{eqnarray*}
&&\|g(x_1),\cdots, g(x_k),g(x_{k+1}),\cdots, g(x_n)\|\\
&&\leq\|g(x_1),\cdots, g(x_k), g(x_{k+1})+\frac{1}{m}g(y_{k+1}),\cdots,g(x_n)+\frac{1}{m}g(y_{n})\|+\frac{1}{m}A_m\\
&&=\frac{1}{m^{n-k}}+\frac{1}{m}A_m,
\end{eqnarray*}
where \begin{eqnarray*}
A_m&=&\sum_{i=1}^{n-k}\|g(x_1),\cdots, g(x_{k+i-1}), g(y_{k+i}),g(x_{k+i+1}),\cdots,g(x_{n})\|+\\
&&\frac{1}{m}\sum_{i=1}^{n-1-k}\|g(x_1),\cdots, g(x_{k+i-1}), g(y_{k+i}),g(y_{k+i+1}),g(x_{k+i+2}),\cdots,g(x_{n})\|+\\
&&\frac{1}{m^{n-k-1}}\|g(x_1),\cdots, g(x_{k}), g(y_{k+1}),g(y_{k+2}),\cdots,g(y_n)\|.
\end{eqnarray*}
Letting ${m\rightarrow+\infty}$ we get the desired equation
\begin{equation*}
\| g(x_1), g(x_2), \cdots, g(x_n)\|=0.
\end{equation*}
\end{prf}

Since it has been showed that if $f$ preserves a fixed nonzero weak $n$-distance and 2-collinearity for the midpoint of a segment then $f(x)-f(0)$ is additive, it is natural to think of such mappings not far from being affine. It is clearly easy to prove $f$ to be an $n$-isometry if it is affine. However it may not be an immediate result since continuity is not implied by preserving nonzero weak $n$-distance.

\begin{prop}\label{isometry}
Let $X$ and $Y$ be two $n$-normed spaces. If $f : X \rightarrow Y$ preserves $w$-$\rho$-$n$-distance for some $\rho>0$ and preserves 2-collinearity for the midpoint of a segment, then $f$ is an affine $n$-isometry.
\end{prop}

\begin{prf}We first prove that $f$ is affine. For this purpose, we only need to show that the mapping $g: X\rightarrow Y$ defined by $g(x)=f(x)-f(0)$ is linear. By Lemmas \ref{add} and \ref{col}, the mapping $g$ is injective, additive and preserves 2-collinearity.
Let $x\in X$ with $x\neq 0$ and $t\in \mathbb{R}$ with $t\neq 1$. Since $0, x, tx$ are collinear, there exists a unique real number $s$ such that $g(tx) = sg(x).$ We can define $\phi: \mathbb{R}\rightarrow \mathbb{R}$ by $\phi(t)=s$ i.e., $$g(tx) = \phi(t) g(x), \ \forall t\in \mathbb{R}.$$ Then clearly, the mapping $\phi$ is injective, additive with $\phi(0)=0$ and $\phi(1)=1$.  Moreover $\phi$ does not depend on the choice of $x$ under the assumption of linear independence. Indeed, choose $y\in X$ such that $x$ and $y$ are linearly independent and let $\phi_1: \mathbb{R}\rightarrow \mathbb{R}$ be a mapping such that $$g(ty) = \phi_1(t) g(y), \ \forall t\in \mathbb{R}.$$ Since $0, x+y, t(x+y)$ are collinear, $$0,  g(x)+g(y), \phi(t)g(x)+\phi_1(t)g(y)$$ are collinear. Note that if
$g(x)$ and $g(y)$ are linearly independent, then $\phi(t)=\phi_1(t)$, as desired. In fact, if $n>2$, there exist $x_3, \cdots, x_n \in X$ such that \begin{eqnarray*}
\| g(x), g(y), g(x_3), \cdots, g(x_n)\| = \|x,y, x_3, \cdots, x_n\|=\rho.\end{eqnarray*} Then $g(x)$ and $g(y)$ are linearly independent. If $n=2$, choose a real number $a$ such that $$\|g(x), g(ay)\|=\|x, ay\|=\rho.$$ Then $g(x)$ and $g(ay)$ are linearly independent, and thus so are $g(x)$ and $g(y)$.
We will prove that $\phi$ is an endomorphism. For any $t, s\in \mathbb{R}$, $0, x+sy, tx+tsy$ are collinear, and then $$0, g(x)+\phi(s)g(y), \phi(t)g(x)+\phi(ts)g(y)$$ are collinear. It follows that $\phi(st)=\phi(s)\phi(t)$ for any $t, s\in \mathbb{R}$. It is well-known that the every nonzero endomorphism of $\mathbb{R}$ is the identity. Then for any $x\in X$ and $t\in \mathbb{R}$, $g(t x) = tg(x)$. Thus $g$ is linear. It is easy to see that $g$ is an $n$-isometry, and hence so is $f$. The proof is complete.
\end{prf}

\begin{rem}\emph{Proposition \ref{isometry} has been shown in \cite[Lemma 3.4]{Ma2}.  Unfortunately the proof given in \cite[Lemma 3.4]{Ma2} contains a mistake. The statement ``$\lim_{k\rightarrow \infty}\|g(rx)-g(r_kx),g(x_2^k),g(x_3^k),\cdots,g(x_n^k)\|=0$ (pp 978, line 11 of \cite{Ma2})'' could not be obtained from the discussing proof in \cite{Ma2}.  For a counterexample, consider $g$ to be the identity, i.e., $g(x)=x$ for every $x\in X$. We may assume that $r$ is an irrational number since the rational case is settled. For each $k$, choose $x_2^k,\cdots,x_n^k$ such that $\|x,x_2^k,x_3^k,\cdots,x_n^k\|=(2+[|r-r_k|])/|r-r_k|$. Then clearly $\|x, x_2^k,x_3^k,\cdots,x_n^k\|>1$ and $|r-r_k|\cdot\|x,x_2^k,x_3^k,\cdots,x_n^k\|=2+[|r-r_k|]$ is a rational number as required in \cite{Ma2}. However, $\|g(rx)-g(r_kx),g(x_2^k),g(x_3^k),\cdots,g(x_n^k)\|=|r-r_k|\cdot\|x,x_2^k,x_3^k,\cdots,x_n^k\|>1$ for every $k$. Therefore the limit cannot be 0 as $k$ goes to infinity. The remaining results Lemma 3.5, Theorem 3.6, Corollary 3.7 and Corollary 3.8 in  \cite{Ma2} following from the the main lemma 3.4 need a new proof. For this and our main result (Theorem 2.6), we hence include a different proof in this paper.}
\end{rem}

We are now ready to prove our main result that gives a positive answer to the Aleksandrov problem in $n$-normed spaces. For a real vector space $X$, we denote the line joining two different points $x,y \in X$ by $\overline{xy}$ and affine$(M)$ by the affine subspace generated by $M\subset X$, respectively.

\begin{thm}
Let $X$ and $Y$ be two $n$-normed spaces. If a surjective mapping $f : X \rightarrow Y$ has $n$-DOPP, then $f$ is an affine $n$-isometry.
\end{thm}
\begin{prf}In the following proof, without loss of generality we can assume that $f(0)=0$. We first prove that $f^{-1}$ preserves 2-collinearity. This is equivalent to showing that if $x, y , z \in X$ are not collinear then $f(x), f(y), f(z)$ are not collinear. Indeed, choose $x_3, \cdots, x_n \in X$ such that $r=\| x-z, y-z, x_3-z, \cdots, x_n-z\|\neq 0$. Set $$u:=z+\frac{(x-z)+(y-z)}{r}.$$ It is easy to check that \begin{eqnarray*}\|x-z, u-z, x_3-z, \cdots, x_n-z\|=\|y-z, u-z, x_3-z, \cdots, x_n-z\|=1.\end{eqnarray*} Since $f$ has $n$-DOPP, \begin{eqnarray}\label{d1}&&\|f(x)-f(z), f(u)-f(z), f(x_3)-f(z), \cdots, f(x_n)-f(z)\|=1\\&&\label{d2}\|f(y)-f(z), f(u)-f(z), f(x_3)-f(z), \cdots, f(x_n)-f(z)\|=1.\end{eqnarray} If there exists some $t\in \mathbb{R}$ such that $f(x)-f(z)=t(f(y)-f(z))$. By (\ref{d1}),(\ref{d2}) and since $f$ is injective,
we obtain that $t=-1$ and so $f(z)=(f(x)+f(y))/2$. Similarly, $f(x)=(f(z)+f(y))/2$. It follows that $f(x)=f(y)=f(z)$, which is impossible.

To see our conclusion, we shall show that $f$ preserves 2-collinearity for the midpoint of a segment. If this does not hold, then there exist $x\neq y \in X$ with $z=(x+y)/2$ such that $f(x), f(y), f(z)$ are not collinear. Now let $w\in X$ such that $$f(w)=\frac{f(x)+f(y)}{2}.$$ Since $f^{-1}$ preserves 2-collinearity, there exists a scalar $t$ such that $y-w=t(x-w)$. We can choose $x_2,\cdots,x_n\in X$ satisfying $\|y-w,x_2,\cdots,x_n\|=1$ and $\overline{0x_2}$ intersects $\overline{xy}$ only in one point denoted by $x_0$. We claim that the $f$-image $f(\overline{0x_2})$ belongs to a line $\overline{0f(x_2)}$ in $Y$. Otherwise, there are $u,v\subset\overline{0x_2}$ such that $f(u),f(v),f(x_0)$ are not collinear.
Set
$$ E:=\mbox{affine}(f(u),f(x_0),f(v)) \quad \mbox{and} \quad F:=\mbox{affine}(f(x),f(x_0),f(y),f(z)).$$
Since $f^{-1}$ preserves 2-collinearity, we have $f^{-1}(E)\subset\overline{0x_2}$ and $f^{-1}(F)\subset\overline{xy}$.
Observe that $f(x_0)\in E\cap F$. Then $E\cap F$ contains infinity points. However,
\begin{equation*}
f^{-1}(E\cap F)\subset f^{-1}(E)\cap f^{-1}(F) \subset\overline{0x_2}\cap\overline{xy}=\{x_0\}.
\end{equation*}
A contradiction since $f$ is injective. By the claim, there are scalars $s_1, s_2$ such that $f(tx_2)=s_1f(x_2)$ and $f(-tx_2)=s_2f(x_2)$.
Since $f$ has $n$-DOPP,  we have
 \begin{eqnarray*}
 &&\|f(y)-f(w),f(x_2),\cdots,f(x_n)\|=\frac{1}{2}\|f(x)-f(y),f(x_2),\cdots,f(x_n)\|=1, \\
 &&\|f(x)-f(w),f(tx_2),\cdots,f(x_n)\|=\frac{1}{2}\|f(x)-f(y),s_1f(x_2),\cdots,f(x_n)\|=1,\\
 &&\|f(x)-f(w),f(-tx_2),\cdots,f(x_n)\|=\frac{1}{2}\|f(x)-f(y),s_2f(x_2),\cdots,f(x_n)\|=1.
 \end{eqnarray*}
It follows that $|s_i|=1$. Since $f$ is injective, the only possibility is that $s_{1}=-1$ and $s_{2}=1$. Thus $t=-1$. Therefore $w=(x+y)/2$. A contradiction guarantees that $f$ preserves 2-collinearity for the midpoint of a segment.  Proposition \ref{isometry} thus completes the proof.
\end{prf}

Next, we shall show that the result of W. Benz holds in $n$-strictly convex spaces.

\begin{defn} An $n$-normed space $X$ is said to be \emph{$n$-strictly convex space} if for any $x_0, x_1, \cdots, x_n\in X$, $x_2, \cdots, x_n \notin span\{x_0, x_1\}$ and $\|x_0+x_1, x_2, \cdots, x_n\|=\|x_0, x_2, \cdots, x_n\|+\|x_1, x_2, \cdots, x_n\|>0$ imply $x_0=tx_1$ for some $t\geq0$.
\end{defn}

\begin{thm}\label{benz}
Let $X$ and $Y$ be two $n$-normed spaces, and let $Y$ be $n$-strictly convex. If $f : X \rightarrow Y$ preserves two $n$-distances $\rho$ and $N\rho$ for some $\rho>0$ and some integer $N>1$, then $f$ is an affine $n$-isometry.
\end{thm}

\begin{prf}
It follows from Proposition \ref{isometry} that we need only prove that $f$ preserves 2-collinearity for the midpoint of a segment.

(a) We first prove that $f$ preserves $2\rho$-$n$-distance. Assume that $N>2$ and $f$ preserves $n$-distances $\rho$ and $N\rho$.
Let $x_1, x_2 \cdots, x_n, y_1, y_2 \cdots, y_n $ be in $X$ such that $$\| x_1-y_1, x_2-y_2,\cdots, x_n-y_n\|=2\rho,$$ and set
 \begin{eqnarray*}
\omega_i=y_1+ i(\frac{x_1-y_1}{2}), \quad \forall i\in \mathbb{N}\cup\{0\}.
\end{eqnarray*}
Then $\omega _{0}=y_1$, $\omega_2=x_1$ and
\begin{eqnarray*}
\omega _{i} - \omega _{i-1}= \frac{x_1-y_1}{2}, \quad \forall i\in \mathbb{N}.
\end{eqnarray*}
It follows that \begin{eqnarray*}
\| \omega _{i} - \omega _{i-1}, x_2-y_2,\cdots, x_n-y_n\|=\rho, \quad \forall i\in \mathbb{N}.
\end{eqnarray*} and $\|\omega _{N} - y_1, x_2-y_2,\cdots, x_n-y_n\|=N\rho$.
Since $f$ preserves $\rho$-$n$-distance, by the triangle inequality, we have
\begin{eqnarray*} &&\| f(x_1)-f(y_1), \cdots, f(x_n)-f(y_n) \| \\ &\leq& \| f(\omega_2)-f(\omega_1),\cdots, f(x_n)-f(y_n)\|+\| f(\omega_1)-f(\omega_0), \cdots, f(x_n)-f(y_n) \| \\ &=& 2\rho  \end{eqnarray*} and similarly,
\begin{eqnarray*} \| f(\omega_N)-f(x_1), \cdots, f(x_n)-f(y_n)\| \leq (N-2)\rho. \end{eqnarray*}
Therefore, \begin{eqnarray*} N\rho&=&\| f(\omega_N)-f(y_1), \cdots, f(x_n)-f(y_n)\| \\
&\leq& \| f(\omega_N)-f(x_1), \cdots, f(x_n)-f(y_n)\|+\| f(x_1)-f(y_1), \cdots, f(x_n)-f(y_n) \|\\
&\leq& (N-2)\rho+ 2\rho=N\rho. \end{eqnarray*}
This implies that \begin{eqnarray*} \| f(x_1)-f(y_1), \cdots, f(x_n)-f(y_n) \|=2\rho. \end{eqnarray*}

(b) Let $z=(x+y)/2$ for distinct $x,y\in X$. Let $g(x)=f(x)-f(0)$. Then $g$ preserves two $n$-distances $\rho$ and $2\rho$.
Thus, there is no loss of generality in assuming that $f(0)=0$.
We shall prove that there exist $x_2, x_3, \cdots, x_n \in X$ such that $$\|y-z, x_2, \cdots, x_n\|=\rho$$ and \begin{eqnarray*}f(x_i) \not \in \mbox{span}\{f(y)-f(z), f(x)-f(z)\} \quad \mbox {for} \ i=2,3,\cdots, n.\end{eqnarray*} Choose $y_2, y_3, \cdots, y_n \in X$ such that $\|y-z, y_2, \cdots, y_n\|=\rho$. We define the set $C_2$ to consist of all elements $\nu$ in $X$ such that $\|y-z, \nu, y_3, \cdots, y_n\|=\rho$, that is \begin{eqnarray*}C_2:=\{\nu\in X: \|y-z, \nu, y_3, \cdots, y_n\|=\rho\}.\end{eqnarray*} We can choose $x_2\in C_2$ such that
\begin{eqnarray*}f(x_2) \not \in \mbox{span}\{f(y)-f(z), f(x)-f(z)\}.\end{eqnarray*} Otherwise, assume that for every $\nu\in C_2$  there exist $\alpha, \beta \in \mathbb{R}$ such that
\begin{eqnarray}\label{a1}f(\nu)=\alpha (f(y)-f(z))+ \beta (f(x)-f(z)).\end{eqnarray} Note that $ \|y-z, \nu, y_3 \cdots, y_n\|=\rho$.
It follows that \begin{eqnarray*} \|x-z, \nu, y_3, \cdots, y_n\|=\rho. \end{eqnarray*}
Then
\begin{eqnarray}\label{a2} &&\|f(y)-f(z), f(\nu), f(y_3), \cdots, f(y_n)\|=\rho, \\ \label{a3}&&\|f(x)-f(z), f(\nu), f(y_3), \cdots, f(y_n)\|=\rho, \end{eqnarray}
It follows from (\ref{a1}), (\ref{a2}) and (\ref{a3}) that
\begin{eqnarray*}
&&|\beta| \| f(y)-f(z), f(x)-f(z),  \cdots, f(y_n)\|=\rho\\
&&|\alpha|\| f(x)-f(z), f(y)-f(z),  \cdots, f(y_n)\|=\rho.
\end{eqnarray*}
This yields $|\alpha|=|\beta|.$ Moreover, $|\alpha|$ is a fixed positive real number. Therefore, there are at most four elements in $f(C_2)$. This is impossible, because the set $C_2$ contains ``enough'' elements.  This follows from Lemma \ref{add} that  $f$ is injective and for each $r\in \mathbb{R}$, the element $\nu_r:=y_2+r(y-z)$ belongs to $C_2$. So there exists $x_2\in C_2$ such that
\begin{eqnarray*}f(x_2) \not \in \mbox{span}\,\{f(y)-f(z), f(x)-f(z)\}.\end{eqnarray*} Next, set \begin{eqnarray*}C_3:=\{\nu\in X: \|y-z, x_2, \nu, y_4, \cdots, y_n\|=\rho\}.\end{eqnarray*} By the same method as above, we can choose $x_3\in C_3$ such that
\begin{eqnarray*}f(x_3) \not \in \mbox {span}\{f(y)-f(z), f(x)-f(z)\}.\end{eqnarray*} This process can be repeated until we obtain the promised $x_2, x_3, \cdots, x_n \in X$ such that $\|y-z, x_2, \cdots, x_n\|=\rho$ and \begin{eqnarray*}f(x_i) \not \in \mbox{span}\{f(y)-f(z), f(x)-f(z)\} \ \mbox{ for} \ \ i=2,3,\cdots, n.\end{eqnarray*}

(c) We are now  ready to show the desired result that $f$ preserves 2-collinearity for the the midpoint of a segment.
Let $z=(x+y)/{2}$ for distinct $x,y\in X$. Let $x_2, x_3, \cdots, x_n $ be in $X$ such that
\begin{eqnarray*}\|y-z, x_2, \cdots, x_n\|=\rho\end{eqnarray*} and
\begin{eqnarray*}f(x_i) \not \in \mbox {span}\{f(y)-f(z), f(x)-f(z)\} \ \mbox{for}  \ i=2,3,\cdots, n.\end{eqnarray*} Then we deduce from the fact that $f$ preserves $n$-distances $\rho$ and $2\rho$ that
\begin{eqnarray*}
&&\| f(y)-f(x), f(x_2),  \cdots, f(x_n)\|\\&=&\|f(y)-f(z),f(x_2),  \cdots, f(x_n)\|+\|f(x)-f(z), f(x_2),  \cdots, f(x_n)\|.
\end{eqnarray*}
Since $Y$ is $n$-strictly convex, there exists a real number $t>0$ such that $$f(y)-f(z)=t(f(z)-f(x)).$$
This completes the proof.
\end{prf}

\begin{rem}
\emph{ \cite[Theorem 11]{Ma} tried to generalize Benz's Theorem on $n$-normed spaces. However, on the part (d) of the proof of \cite[Theorem 11]{Ma}
the statement that $f(p_2)-f(p_1)=t(f(p_1)-f(p_0))$ for some $t$ cannot follow just from
\begin{align*}
& \| f(p_2)-f(p_0), f(y_2)-f(x), \cdots,f(y_n)-f(x)\|\\
&=\| f(p_2)-f(p_1), f(y_2)-f(x), \cdots, f(y_n)-f(x)\|
\\&+\| f(p_1)-f(p_0), f(y_2)-f(x), \cdots, f(y_n)-f(x)\| = 2\rho.
\end{align*}
It remains to check that $f(y_i)-f(x)\not \in \mbox{span}\{f(p_2)-f(p_0), f(p_1)-f(p_0)\} \quad \mbox {for} \ i=2,3,\cdots, n$ (It is the demand from the definition of $n$ strictly convexity (\cite[definition 3]{Ma} or Definition 2.7 of our paper)). It is a hard and key step which cannot be missed.  }
\end{rem}

\subsection*{Acknowledgements}
The authors wish to express their appreciation to Guanggui Ding for many very helpful comments regarding isometric theory in Banach spaces.


\begin{thebibliography}{1000}

\bibitem{A} A. D. Alekandrov, Mappings of families of sets, {\it Soviet Math. Dokl.} 11 (1970), 116--120.

\bibitem{B} W. Benz, Isometrien in normierten R\"{a}umen, {\it Aequationes Math.}29 (1985), 204--209.

\bibitem{BB} W. Benz and H.Berens, A contribution to a theorem of Ulam and Mazur, {\it Aequationes Math.}  34 (1987), 61--63.

\bibitem{C} H. Y. Chu, On the Mazur--Ulam problem in linear 2-normed spaces,  {\it J. Math. Anal. Appl.} 327 (2007), 1041--1045.

\bibitem{HPP} H. Y. Chu, C. G. Park and W. G. Park, The Aleksandrov problem in linear 2-normed spaces, {\it J. Math. Anal. Appl.} 289 (2004), 666--672.

\bibitem{CKK} H. Y. Chu, S. K. Choi and D. S. Kang, Mappings of conservative distances in linear $n$-normed spaces, {\it Nonlinear Anal.} 70 (2009), 1168--1174.

\bibitem{CLP} H. Chu, K. Lee and C. Park, On the Aleksandrov problem in linear $n$-normed spaces, {\it Nonlinear Anal.} 59 (2004), 1001--1011.

\bibitem{G1} S. G\"{a}hler, Lineare 2-normierte R\"{a}ume, {\it  Math. Nachr. } 28 (1964), 1--43.

\bibitem{G} S. G\"{a}hler, Untersuchungen \"{u}ber verallgemeinerte m-metrische R\"{a}ume, {\it I. Math. Nachr.} 40 (1969), 165--189.

\bibitem{G2}Gy. P. Geh\'{e}r, A contribution to the Aleksandrov conservative distance problem in two dimensions, {\it Linear Algebra Appl.}, 481(2015), 280--287.

\bibitem{GH} H. Gunawan, The space of $p$-summable sequences and its natural $n$-norm. {\it Bull. Austral. Math. Soc.} 64 (2001), 137--147.

\bibitem{Ma}  Y. Ma, The Aleksandrov-Ben-Rassias problem on linear $n$-normed spaces. Monatshefte f\"{u}r Mathematik, 2015, 1-12.

\bibitem{Ma2} Y. Ma, Isometry on linear $n$-normed spaces, {\it Ann. Acad. Sci. Fenn. Math.} 39 (2) (2014), 973--981.

\bibitem{M1} A. Misiak, $n$-inner product spaces, {\it Math. Nachr.} 140 (1989), 299--319.

\bibitem{M2} A. Misiak, Orthogonality and orthogonormality in $n$-inner product spaces, {\it Math. Nachr.} 143 (1989), 249--261.

\bibitem{R1} Th. M. Rassias, Mappings that preserve unit distance, {\it Indian J. Math.} 32 (1990), 275--278.

\bibitem{R2} Th. M. Rassias, Properties of isometric mappings, {\it J. Math. Anal. Appl.} 235 (1997), 108--121.

\bibitem{R3} Th. M. Rassias, On the  Aleksandrov problem of conservative distances and the Mazur-Ulam theorem, {\it Nonlinear Anal.} 47 (2001), 2597--2608.

\bibitem{R4} Th. M. Rassias, On the Aleksandrov problem for isometric mappings, {\it Appl. Anal. Discrete Math.} 1 (2007),  18--28.

\bibitem{RS} Th. M. Rassias and P. \v{S}emrl, On the Mazur-Ulam problem and the Aleksandrov problem for unit distance preserving mappings, {\it Proc. Amer. Math. Soc. }118 (1993), 919--925.


\end{thebibliography}
\end{document}